\documentclass[12pt]{article}%
\usepackage{amsfonts}
\usepackage{sw20bams}
\usepackage{amsmath}
\usepackage{amssymb}
\usepackage{graphicx}%
\providecommand{\U}[1]{\protect\rule{.1in}{.1in}}
\begin{document}

\title{Random Cyclic Quadrilaterals}
\author{Steven Finch}
\date{October 3, 2016}
\maketitle

\begin{abstract}
The circumcircle of a planar convex polygon $P$ is a circle $C$ that passes
through all vertices of $P$. If such a $C$ exists, then $P$ is said to be
cyclic. Fix $C$ to have unit radius. While any two angles of a uniform cyclic
triangle are negatively correlated, any two sides are independent. In
contrast, for a uniform cyclic quadrilateral, any two sides are negatively
correlated, whereas any two adjacent angles are uncorrelated yet dependent.

\end{abstract}

\footnotetext{Copyright \copyright \ 2016 by Steven R. Finch. All rights
reserved.}To generate a cyclic triangle is easy: select three independent
uniform points on the unit circle and connect them. \ To generate a cyclic
quadrilateral is harder: select four such points and connect them in, say, a
counterclockwise manner. \ Convexity follows immediately \cite{Pi-quadrilat},
as does the fact that opposite angles are supplementary \cite{M1-quadrilat,
M2-quadrilat}. The inter-relationship of adjacent angles is more mysterious,
as we shall soon see. \ Our initial focus, however, will be on adjacent sides,
opposite sides and diagonals.

Let the four vertices be given by $\exp(i\,\theta_{k})$, where $i$ is the
imaginary unit, $0\leq\theta_{1}<$ $\theta_{2}<\theta_{3}<\theta_{4}<2\pi$ are
central angles relative to the horizontal axis, and $1\leq k\leq4$. \ Define
$\theta_{0}=\theta_{4}-2\pi$ and $\theta_{5}=\theta_{1}+2\pi$ for convenience,
then polygonal sides $s_{k}$ and polygonal angles $\alpha_{k}$ are given by%
\[%
\begin{array}
[c]{ccc}%
s_{k}=2\sin\left(  \dfrac{\theta_{k}-\theta_{k-1}}{2}\right)  , &  &
\alpha_{k}=\dfrac{\theta_{k+1}-\theta_{k-1}}{2}.
\end{array}
\]
Proof of the $s_{k}$ expression comes from the Law of Cosines and a half angle
formula:%
\begin{align*}
s_{k}^{2}  &  =1+1-2\cdot1\cdot1\cos(\theta_{k}-\theta_{k-1})=2\left[
1-\cos(\theta_{k}-\theta_{k-1})\right] \\
&  =4\,\frac{1-\cos(\theta_{k}-\theta_{k-1})}{2}=4\sin^{2}\left(
\dfrac{\theta_{k}-\theta_{k-1}}{2}\right)  .
\end{align*}
Proof of the $\alpha_{k}$ expression follows the fact that an inscribed angle
is one-half the length of its intercepted circular arc. \ The polygonal
diagonals $d_{k}$ clearly satisfy $d_{k}=2\sin(\alpha_{k})$. \ Let also
$\omega$ denote the smaller of the two angles at the intersection point
between the diagonals.%
%TCIMACRO{\FRAME{ftbpFU}{5.1647in}{4.5368in}{0pt}{\Qcb{Vertices, angles and
%sides of a cyclic quadrilateral.}}{}{quadrilat.eps}%
%{\special{ language "Scientific Word";  type "GRAPHIC";
%maintain-aspect-ratio TRUE;  display "USEDEF";  valid_file "F";
%width 5.1647in;  height 4.5368in;  depth 0pt;  original-width 12.1558in;
%original-height 10.6614in;  cropleft "0";  croptop "1";  cropright "1";
%cropbottom "0";  filename 'quadrilat.eps';file-properties "XNPEU";}} }%
%BeginExpansion
\begin{figure}[ptb]%
\centering
\includegraphics[
height=4.5368in,
width=5.1647in
]%
{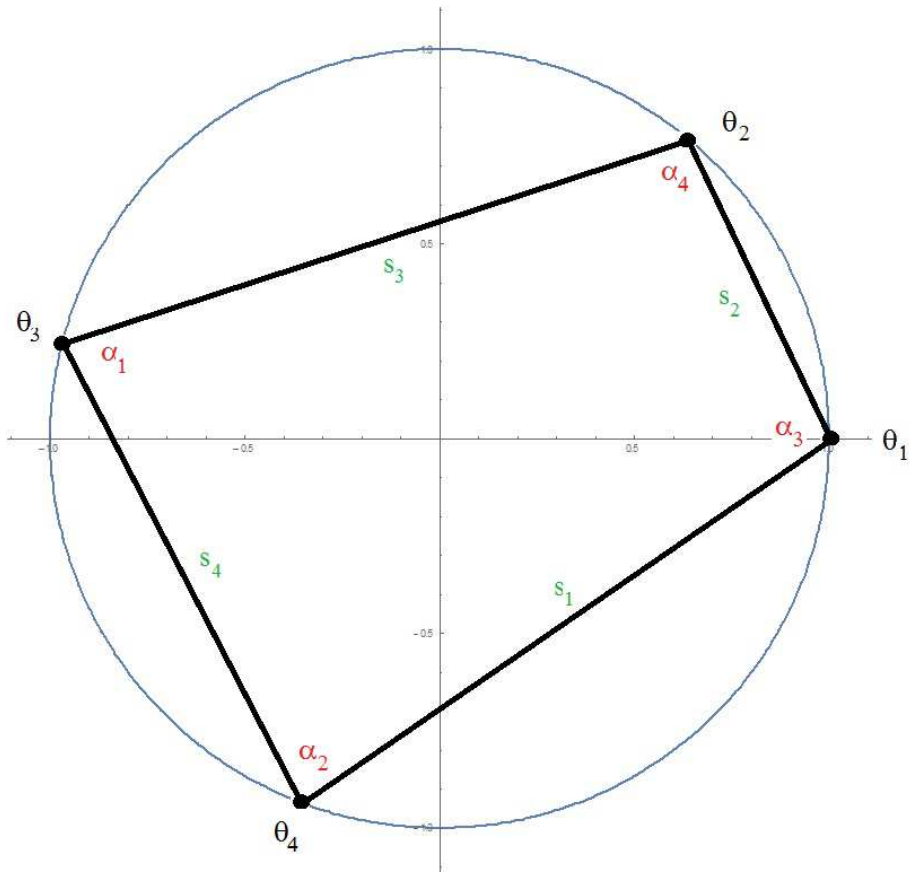}%
\caption{Vertices, angles and sides of a cyclic quadrilateral.}%
\end{figure}
%EndExpansion

Our labor draws upon the distribution of the order statistics $\theta_{1}$,
$\theta_{2}$, $\theta_{3}$, $\theta_{4}$. \ We must be careful in summarizing
the results because, while $s_{2}$, $s_{3}$, $s_{4}$ possess the same
distribution, the one corresponding to $s_{1}$ is different. \ Hence, to make
statements regarding arbitrary sides $s$, $t$, $u$, $v$ of the quadrilateral,
we must use a $(3/4,1/4)$-mixture of densities. \ Likewise, $\alpha_{2}$ and
$\alpha_{1}$ possess distinct distributions. \ Thus, to make statements
regarding arbitrary adjacent angles $\alpha$, $\beta$ of the quadrilateral, we
must use a $(1/2,1/2)$-mixture of densities.

A probabilistic analysis of the perimeter $s+t+u+v$ and area $2\sin
(\alpha)\sin(\beta)\sin(\omega)$ is beyond our current capabilities.
\ Hopefully the groundwork established here will be a launching point for
someone else's research in the near future.

\section{Sides}

Let $X_{1}<X_{2}<X_{3}<X_{4}$ denote the order statistics for a random sample
of size $4$ from the uniform distribution on $[0,1]$. \ The density for
$(X_{1},X_{2})=(x,y)$ is \cite{Gi-quadrilat, HN-quadrilat}%
\[
\left\{
\begin{array}
[c]{lll}%
12(1-y)^{2} &  & \text{if }0<x<y<1,\\
0 &  & \text{otherwise;}%
\end{array}
\right.
\]
the density for $(X_{1},X_{3})=(x,y)$ is%
\[
\left\{
\begin{array}
[c]{lll}%
24(y-x)(1-y) &  & \text{if }0<x<y<1,\\
0 &  & \text{otherwise;}%
\end{array}
\right.
\]
the density for $(X_{1},X_{4})=(x,y)$ is%
\[
\left\{
\begin{array}
[c]{lll}%
12(y-x)^{2} &  & \text{if }0<x<y<1,\\
0 &  & \text{otherwise;}%
\end{array}
\right.
\]
the density for $(X_{2},X_{4})=(x,y)$ is%
\[
\left\{
\begin{array}
[c]{lll}%
24x(y-x) &  & \text{if }0<x<y<1,\\
0 &  & \text{otherwise.}%
\end{array}
\right.
\]
Consider the transformation $(x,y)\mapsto(y-x,y)=(u,v)$. \ Since this has
Jacobian determinant $1$ and since $0<u<v<1$, it follows that the density for
$X_{2}-X_{1}$ is%
\[
12%
%TCIMACRO{\dint \limits_{u}^{1}}%
%BeginExpansion
{\displaystyle\int\limits_{u}^{1}}
%EndExpansion
(1-v)^{2}dv=\left.  -4(1-v)^{3}\right\vert _{u}^{1}=4(1-u)^{3};
\]
the density for $X_{3}-X_{1}$ is%
\[
24%
%TCIMACRO{\dint \limits_{u}^{1}}%
%BeginExpansion
{\displaystyle\int\limits_{u}^{1}}
%EndExpansion
u(1-v)dv=\left.  -12u(1-v)^{2}\right\vert _{u}^{1}=12u(1-u)^{2};
\]
the density for $X_{4}-X_{1}$ is%
\[
12%
%TCIMACRO{\dint \limits_{u}^{1}}%
%BeginExpansion
{\displaystyle\int\limits_{u}^{1}}
%EndExpansion
u^{2}dv=\left.  -12u^{2}(1-v)\right\vert _{u}^{1}=12u^{2}(1-u);
\]
the density for $X_{4}-X_{2}$ is%
\[
24%
%TCIMACRO{\dint \limits_{u}^{1}}%
%BeginExpansion
{\displaystyle\int\limits_{u}^{1}}
%EndExpansion
u(v-u)dv=\left.  12u(v-u)^{2}\right\vert _{u}^{1}=12u(1-u)^{2}.
\]
We disregard $X_{4}-X_{2}$ further since its distribution is the same as that
for $X_{3}-X_{1}$. \ Consider the scaling $u\mapsto\pi\,u=x$. \ It follows
that%
\[%
\begin{array}
[c]{ccc}%
\text{the density for }\dfrac{\theta_{2}-\theta_{1}}{2}\text{ is} &  &
\dfrac{4}{\pi}\left(  1-\dfrac{x}{\pi}\right)  ^{3};\\
\text{the density for }\dfrac{\theta_{3}-\theta_{1}}{2}\text{ is} &  &
\dfrac{12}{\pi}\left(  \dfrac{x}{\pi}\right)  \left(  1-\dfrac{x}{\pi}\right)
^{2};\\
\text{the density for }\dfrac{\theta_{4}-\theta_{1}}{2}\text{ is} &  &
\dfrac{12}{\pi}\left(  \dfrac{x}{\pi}\right)  ^{2}\left(  1-\dfrac{x}{\pi
}\right)  .
\end{array}
\]
Next, the function $x\mapsto\sin(x)=y$ possesses two preimages $\arcsin(y) $
and $\pi-\arcsin(y)$ in the interval $[0,\pi]$ and has derivative
$\cos(x)=\sqrt{1-y^{2}}$. It follows that the three densities are
\cite{Pa-quadrilat}
\[
\dfrac{4}{\pi\sqrt{1-y^{2}}}\left[  \left(  1-\dfrac{\arcsin(y)}{\pi}\right)
^{3}+\left(  \dfrac{\arcsin(y)}{\pi}\right)  ^{3}\right]  ,
\]%
\[
\dfrac{12}{\pi\sqrt{1-y^{2}}}\left[  \left(  \dfrac{\arcsin(y)}{\pi}\right)
\left(  1-\dfrac{\arcsin(y)}{\pi}\right)  ^{2}+\left(  1-\dfrac{\arcsin
(y)}{\pi}\right)  \left(  \dfrac{\arcsin(y)}{\pi}\right)  ^{2}\right]  ,
\]%
\[
\dfrac{12}{\pi\sqrt{1-y^{2}}}\left[  \left(  \dfrac{\arcsin(y)}{\pi}\right)
^{2}\left(  1-\dfrac{\arcsin(y)}{\pi}\right)  +\left(  1-\dfrac{\arcsin
(y)}{\pi}\right)  ^{2}\left(  \dfrac{\arcsin(y)}{\pi}\right)  \right]
\]
respectively. \ The second and third expressions are identical. \ Finally, the
scaling $y\mapsto2\,y=z$ and an algebraic expansion gives the density for
$s_{2}$ as%
\[
\dfrac{4}{\pi\sqrt{4-z^{2}}}\left[  1-\frac{3\arcsin(\frac{z}{2})\left(
\pi-\arcsin(\frac{z}{2})\right)  }{\pi^{2}}\right]
\]
and the density for both $d_{2}$ and $s_{1}$ as%
\[
\dfrac{4}{\pi\sqrt{4-z^{2}}}\left[  0+\frac{3\arcsin(\frac{z}{2})\left(
\pi-\arcsin(\frac{z}{2})\right)  }{\pi^{2}}\right]  .
\]
We omit details for $s_{3}$, $s_{4}$ (same as $s_{2}$) and $d_{1}$ (same as
$d_{2}$).

Mixing the densities for $s_{2}$ (with weight $3/4$) and for $s_{1}$ (with
weight $1/4$), the density for an arbitrary side $0\leq s\leq2$ emerges:%
\[
\dfrac{3}{\pi\sqrt{4-s^{2}}}\left[  1-\frac{2\arcsin(\frac{s}{2})\left(
\pi-\arcsin(\frac{s}{2})\right)  }{\pi^{2}}\right]
\]
which implies that%
\[%
\begin{array}
[c]{ccc}%
\operatorname*{E}\left(  s\right)  =\dfrac{6}{\pi}-\dfrac{24}{\pi^{3}}, &  &
\operatorname*{E}\left(  s^{2}\right)  =2-\dfrac{3}{\pi^{2}}.
\end{array}
\]
The corresponding moments for a diagonal $0\leq d\leq2$ are $48/\pi^{3}$ and
$2+6/\pi^{2}$. \ Joint moments are available via the joint density of
$\theta_{1}$, $\theta_{2}$, $\theta_{3}$, $\theta_{4}$:%
\[
\left\{
\begin{array}
[c]{lll}%
\dfrac{4!}{(2\pi)^{4}} &  & \text{if }0\leq\theta_{1}<\theta_{2}<\theta
_{3}<\theta_{4}<2\pi,\\
0 &  & \text{otherwise.}%
\end{array}
\right.
\]
For example,%
\begin{align*}
\operatorname*{E}\left(  s_{2}s_{3}\right)   &  =\frac{3}{2\pi^{4}}%
%TCIMACRO{\dint \limits_{0}^{2\pi}}%
%BeginExpansion
{\displaystyle\int\limits_{0}^{2\pi}}
%EndExpansion%
%TCIMACRO{\dint \limits_{\theta_{1}}^{2\pi}}%
%BeginExpansion
{\displaystyle\int\limits_{\theta_{1}}^{2\pi}}
%EndExpansion%
%TCIMACRO{\dint \limits_{\theta_{2}}^{2\pi}}%
%BeginExpansion
{\displaystyle\int\limits_{\theta_{2}}^{2\pi}}
%EndExpansion%
%TCIMACRO{\dint \limits_{\theta_{3}}^{2\pi}}%
%BeginExpansion
{\displaystyle\int\limits_{\theta_{3}}^{2\pi}}
%EndExpansion
\left[  2\sin\left(  \dfrac{\theta_{2}-\theta_{1}}{2}\right)  \right]  \left[
2\sin\left(  \dfrac{\theta_{3}-\theta_{2}}{2}\right)  \right]  d\theta
_{4}d\theta_{3}d\theta_{2}d\theta_{1}\\
&  =\dfrac{48}{\pi^{2}}-\dfrac{384}{\pi^{4}}%
\end{align*}
(same for $\operatorname*{E}\left(  s_{3}s_{4}\right)  $) and
\begin{align*}
\operatorname*{E}\left(  s_{1}s_{2}\right)   &  =\frac{3}{2\pi^{4}}%
%TCIMACRO{\dint \limits_{0}^{2\pi}}%
%BeginExpansion
{\displaystyle\int\limits_{0}^{2\pi}}
%EndExpansion%
%TCIMACRO{\dint \limits_{\theta_{1}}^{2\pi}}%
%BeginExpansion
{\displaystyle\int\limits_{\theta_{1}}^{2\pi}}
%EndExpansion%
%TCIMACRO{\dint \limits_{\theta_{2}}^{2\pi}}%
%BeginExpansion
{\displaystyle\int\limits_{\theta_{2}}^{2\pi}}
%EndExpansion%
%TCIMACRO{\dint \limits_{\theta_{3}}^{2\pi}}%
%BeginExpansion
{\displaystyle\int\limits_{\theta_{3}}^{2\pi}}
%EndExpansion
\left[  2\sin\left(  \dfrac{\theta_{4}-\theta_{1}}{2}\right)  \right]  \left[
2\sin\left(  \dfrac{\theta_{2}-\theta_{1}}{2}\right)  \right]  d\theta
_{4}d\theta_{3}d\theta_{2}d\theta_{1}\\
&  =-\dfrac{24}{\pi^{2}}+\dfrac{384}{\pi^{4}}%
\end{align*}
(same for $\operatorname*{E}\left(  s_{4}s_{1}\right)  $) imply that, for
arbitrary adjacent sides $s$ and $t$,%
\[%
\begin{array}
[c]{ccc}%
\operatorname*{E}\left(  s\,t\right)  =\dfrac{12}{\pi^{2}}, &  &
\rho(s,t)\approx-0.183.
\end{array}
\]
We used
\[
\dfrac{\theta_{1}-\theta_{0}}{2}=\dfrac{\theta_{1}-\left(  \theta_{4}%
-2\pi\right)  }{2}=\dfrac{2\pi-\left(  \theta_{4}-\theta_{1}\right)  }{2}%
=\pi-\dfrac{\theta_{4}-\theta_{1}}{2}%
\]
and $\sin(\pi-z)=\sin(z)$ in writing the preceding integral. \ The same value
$12/\pi^{2}$ is also obtained for the expected product of arbitrary opposite
sides $s$ and $t$. \ The proximity of quadrilateral sides is (evidently)
immaterial when assessing their correlation.%
%TCIMACRO{\FRAME{ftbpFU}{2.919in}{2.9149in}{0pt}{\Qcb{Density function for side
%$s$ in Section 1.}}{}{sides.eps}{\special{ language "Scientific Word";
%type "GRAPHIC";  maintain-aspect-ratio TRUE;  display "USEDEF";
%valid_file "F";  width 2.919in;  height 2.9149in;  depth 0pt;
%original-width 11.1996in;  original-height 11.1855in;  cropleft "0";
%croptop "1";  cropright "1";  cropbottom "0";
%filename 'sides.eps';file-properties "XNPEU";}} }%
%BeginExpansion
\begin{figure}[ptb]%
\centering
\includegraphics[
height=2.9149in,
width=2.919in
]%
{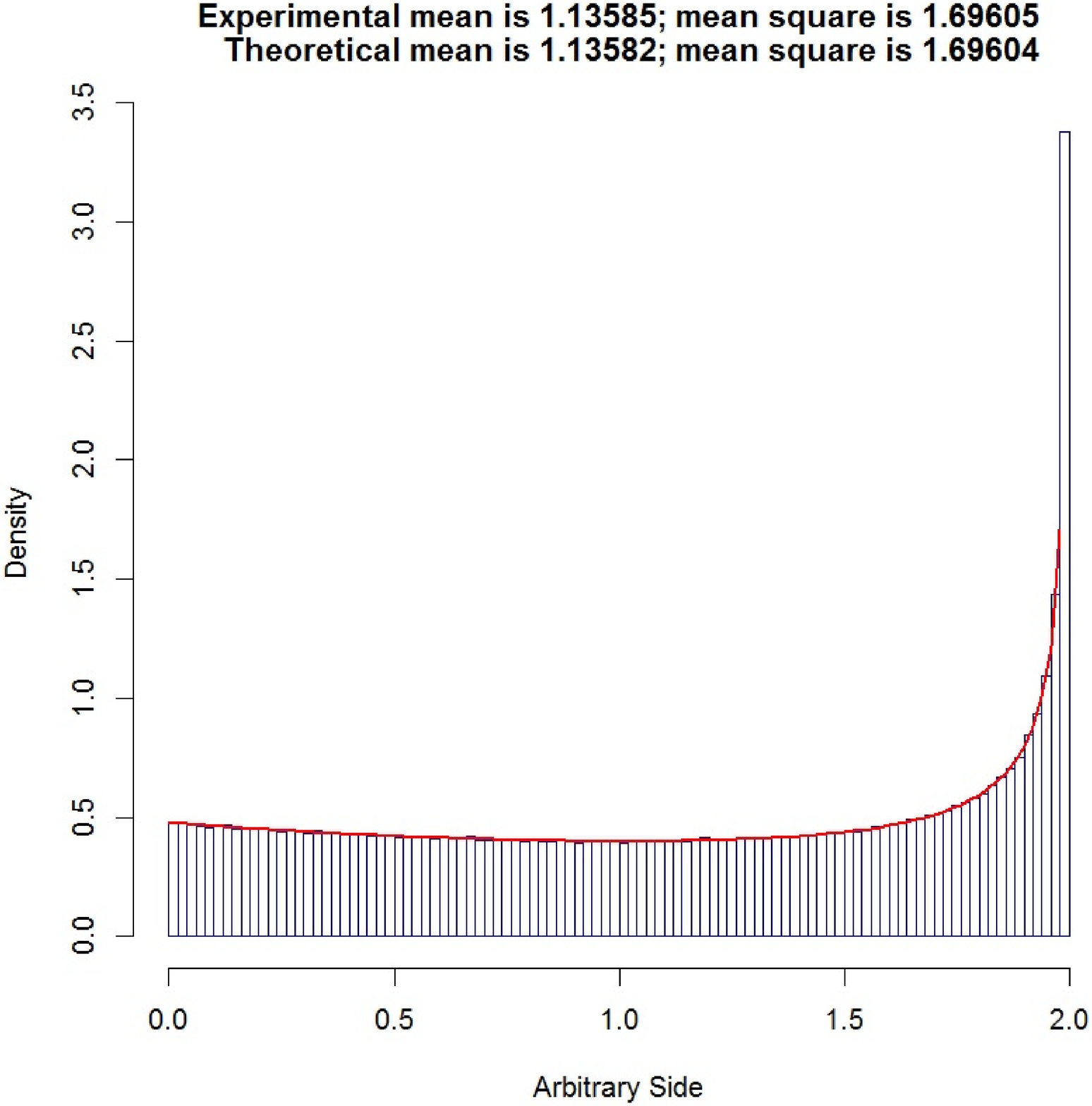}%
\caption{Density function for side $s$ in Section 1.}%
\end{figure}
%EndExpansion%
%TCIMACRO{\FRAME{ftbpFU}{2.919in}{2.9149in}{0pt}{\Qcb{Density function for
%diagonal $d$ in Section 1.}}{}{diagonals.eps}%
%{\special{ language "Scientific Word";  type "GRAPHIC";
%maintain-aspect-ratio TRUE;  display "USEDEF";  valid_file "F";
%width 2.919in;  height 2.9149in;  depth 0pt;  original-width 11.1996in;
%original-height 11.1855in;  cropleft "0";  croptop "1";  cropright "1";
%cropbottom "0";  filename 'diagonals.eps';file-properties "XNPEU";}} }%
%BeginExpansion
\begin{figure}[ptb]%
\centering
\includegraphics[
height=2.9149in,
width=2.919in
]%
{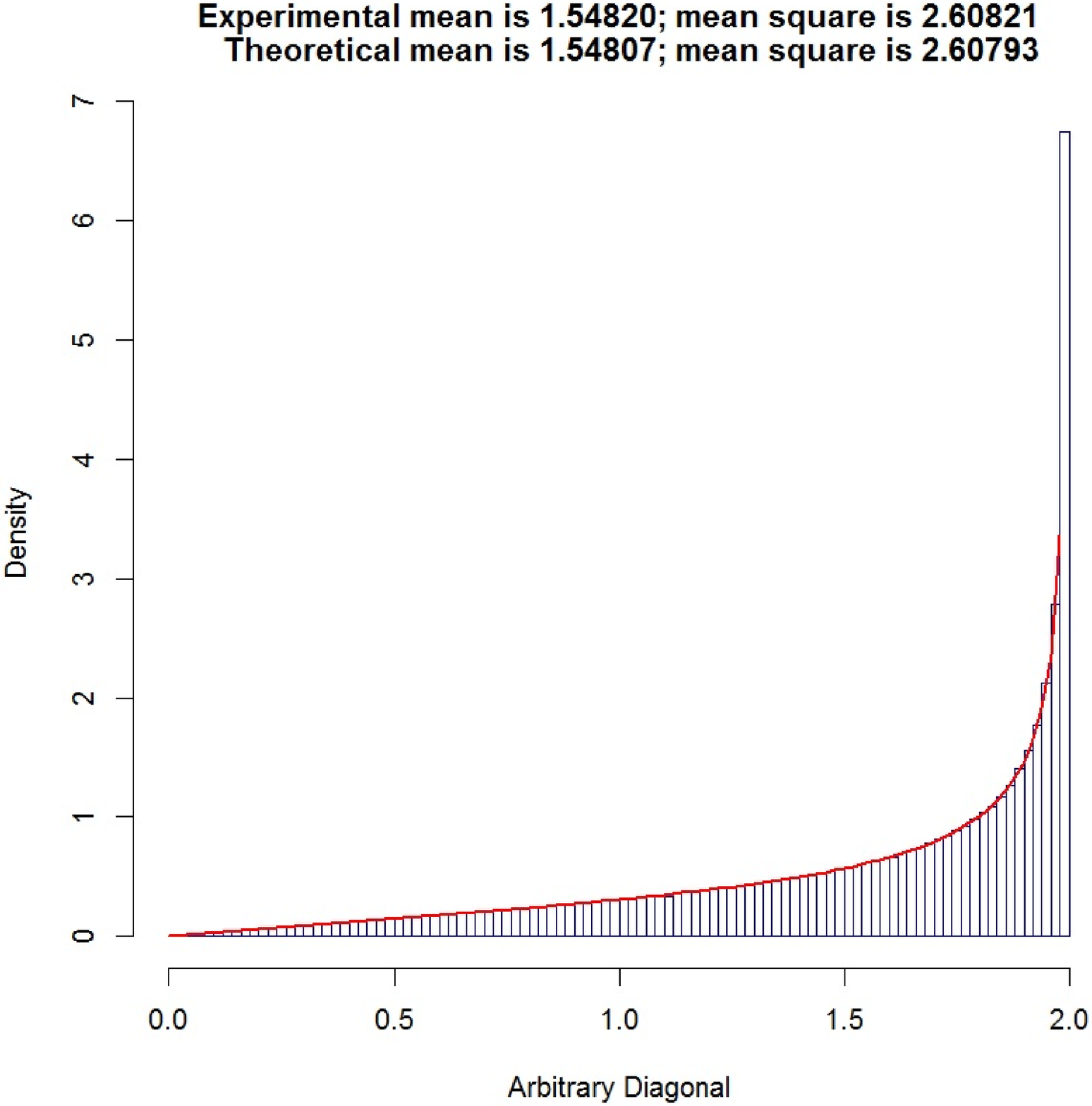}%
\caption{Density function for diagonal $d$ in Section 1.}%
\end{figure}
%EndExpansion

\section{Angles}

By our work starting with $X_{3}-X_{1}$ and $X_{4}-X_{2}$, it is clear that
$\alpha_{2}$ and $\alpha_{3}$ are identically distributed. \ Since $\alpha
_{1}=\pi-\alpha_{3}$, the density of $\alpha_{2}$ is $12x(\pi-x)^{2}/\pi^{4}$
while the density of $\alpha_{1}$ is $12x^{2}(\pi-x)/\pi^{4}$. \ Mixing the
densities for $\alpha_{2}$ and for $\alpha_{1}$ with equal weighting, the
marginal density for an arbitrary angle $0\leq\alpha\leq\pi$ becomes
$6x(\pi-x)/\pi^{3}$.

We need, however, to find the joint distribution for arbitrary adjacent angles
$\alpha$ and $\beta$. \ A fresh approach for obtaining this involves the
Dirichlet$(1,1,1;1)$ distribution on a $3$-dimensional simplex
\cite{Ra-quadrilat, PC-quadrilat, NT-quadrilat}:%
\[
\left\{
\begin{array}
[c]{lll}%
6 &  & \text{if }0<\xi_{1}<1,\text{ }0<\xi_{2}<1,\text{ }0<\xi_{3}<1\text{ and
}\xi_{1}+\xi_{2}+\xi_{3}<1,\\
0 &  & \text{otherwise}%
\end{array}
\right.
\]
and calculation of the joint density for $\eta_{1}=\xi_{1}+\xi_{2}$, $\eta
_{2}=\xi_{2}+\xi_{3}$. \ The list $\xi_{1}$, $\xi_{2}$, $\xi_{3}$ \ can be
thought of as duplicating any one of the eight lists given in Table 1, each
weighted with probability $1/8$. In words, up to the preservation of adjacency
of angles $\pi\,\eta_{1}$, $\pi\,\eta_{2}$, any implicit ordering within
$\xi_{1}$, $\xi_{2}$, $\xi_{3}$ \ has been removed. This formulation will
simplify our work, removing the need to mix distributions (like before) as a
concluding step.

Table 1. \textit{Eight possibilities for }$\xi_{1}$, $\xi_{2}$, $\xi_{3}%
$\textit{.}%
\[%
\begin{tabular}
[c]{|c|c|}\hline
Candidate Lists & Resulting Angles\\\hline
$%
\begin{array}
[c]{ccccc}%
\dfrac{\theta_{2}-\theta_{1}}{2\pi}, &  & \dfrac{\theta_{3}-\theta_{2}}{2\pi
}, &  & \dfrac{\theta_{4}-\theta_{3}}{2\pi}%
\end{array}
$ & $%
\begin{array}
[c]{ccc}%
\pi\,\eta_{1}=\alpha_{2}, &  & \pi\,\eta_{2}=\alpha_{3}%
\end{array}
$\\\hline
$%
\begin{array}
[c]{ccccc}%
\dfrac{\theta_{4}-\theta_{3}}{2\pi}, &  & \dfrac{\theta_{3}-\theta_{2}}{2\pi
}, &  & \dfrac{\theta_{2}-\theta_{1}}{2\pi}%
\end{array}
$ & $%
\begin{array}
[c]{ccc}%
\pi\,\eta_{1}=\alpha_{3}, &  & \pi\,\eta_{2}=\alpha_{2}%
\end{array}
$\\\hline
$%
\begin{array}
[c]{ccccc}%
\dfrac{\theta_{3}-\theta_{2}}{2\pi}, &  & \dfrac{\theta_{4}-\theta_{3}}{2\pi
}, &  & \dfrac{\theta_{5}-\theta_{4}}{2\pi}%
\end{array}
$ & $%
\begin{array}
[c]{ccc}%
\pi\,\eta_{1}=\alpha_{3}, &  & \pi\,\eta_{2}=\alpha_{4}%
\end{array}
$\\\hline
$%
\begin{array}
[c]{ccccc}%
\dfrac{\theta_{5}-\theta_{4}}{2\pi}, &  & \dfrac{\theta_{4}-\theta_{3}}{2\pi
}, &  & \dfrac{\theta_{3}-\theta_{2}}{2\pi}%
\end{array}
$ & $%
\begin{array}
[c]{ccc}%
\pi\,\eta_{1}=\alpha_{4}, &  & \pi\,\eta_{2}=\alpha_{3}%
\end{array}
$\\\hline
$%
\begin{array}
[c]{ccccc}%
\dfrac{\theta_{4}-\theta_{3}}{2\pi}, &  & \dfrac{\theta_{5}-\theta_{4}}{2\pi
}, &  & \dfrac{\theta_{2}-\theta_{1}}{2\pi}%
\end{array}
$ & $%
\begin{array}
[c]{ccc}%
\pi\,\eta_{1}=\alpha_{4}, &  & \pi\,\eta_{2}=\alpha_{1}%
\end{array}
$\\\hline
$%
\begin{array}
[c]{ccccc}%
\dfrac{\theta_{2}-\theta_{1}}{2\pi}, &  & \dfrac{\theta_{5}-\theta_{4}}{2\pi
}, &  & \dfrac{\theta_{4}-\theta_{3}}{2\pi}%
\end{array}
$ & $%
\begin{array}
[c]{ccc}%
\pi\,\eta_{1}=\alpha_{1}, &  & \pi\,\eta_{2}=\alpha_{4}%
\end{array}
$\\\hline
$%
\begin{array}
[c]{ccccc}%
\dfrac{\theta_{5}-\theta_{4}}{2\pi}, &  & \dfrac{\theta_{2}-\theta_{1}}{2\pi
}, &  & \dfrac{\theta_{3}-\theta_{2}}{2\pi}%
\end{array}
$ & $%
\begin{array}
[c]{ccc}%
\pi\,\eta_{1}=\alpha_{1}, &  & \pi\,\eta_{2}=\alpha_{2}%
\end{array}
$\\\hline
$%
\begin{array}
[c]{ccccc}%
\dfrac{\theta_{3}-\theta_{2}}{2\pi}, &  & \dfrac{\theta_{2}-\theta_{1}}{2\pi
}, &  & \dfrac{\theta_{5}-\theta_{4}}{2\pi}%
\end{array}
$ & $%
\begin{array}
[c]{ccc}%
\pi\,\eta_{1}=\alpha_{2}, &  & \pi\,\eta_{2}=\alpha_{1}%
\end{array}
$\\\hline
\end{tabular}
\]

Introducing $\eta_{3}=\xi_{3}$, we have%
\[%
\begin{array}
[c]{l}%
\xi_{1}=\eta_{1}-\eta_{2}+\eta_{3},\\
\xi_{2}=\eta_{2}-\eta_{3},\\
\xi_{3}=\eta_{3}%
\end{array}
\]
and calculate the Jacobian determinant to be equal to $1$. \ From%
\[%
\begin{array}
[c]{l}%
0<\eta_{1}-\eta_{2}+\eta_{3}<1,\\
0<\eta_{2}-\eta_{3}<1,\\
0<\eta_{3}<1,\\
0<\eta_{1}+\eta_{3}<1
\end{array}
\]
it follows that%
\[%
\begin{array}
[c]{l}%
-\eta_{1}+\eta_{2}<\eta_{3}<1-\eta_{1}+\eta_{2},\\
-1+\eta_{2}<\eta_{3}<\eta_{2},\\
0<\eta_{3}<1,\\
-\eta_{1}<\eta_{3}<1-\eta_{1}%
\end{array}
\]
hence $\max\{-\eta_{1}+\eta_{2},0\}<\eta_{3}<\min\{\eta_{2},1-\eta_{1}\}$.
\ There are four cases:

\begin{enumerate}
\item[(1.)] If $1-\eta_{2}<\eta_{1}<\eta_{2}$, then $-\eta_{1}+\eta_{2}%
<\eta_{3}<1-\eta_{1}$

\item[(2.)] If $\eta_{1}<\eta_{2}<1-\eta_{1}$, then $-\eta_{1}+\eta_{2}%
<\eta_{3}<\eta_{2}$

\item[(3.)] If $1-\eta_{1}<\eta_{2}<\eta_{1}$, then $0<\eta_{3}<1-\eta_{1}$

\item[(4.)] If $\eta_{2}<\eta_{1}<1-\eta_{2}$, then $0<\eta_{3}<\eta_{2} $
\end{enumerate}

\noindent giving rise to%
\[%
\begin{array}
[c]{cc}%
%TCIMACRO{\dint \limits_{-\eta_{1}+\eta_{2}}^{1-\eta_{1}}}%
%BeginExpansion
{\displaystyle\int\limits_{-\eta_{1}+\eta_{2}}^{1-\eta_{1}}}
%EndExpansion
6\,d\eta_{3}=6\left(  1-\eta_{2}\right)  , &
%TCIMACRO{\dint \limits_{-\eta_{1}+\eta_{2}}^{\eta_{2}}}%
%BeginExpansion
{\displaystyle\int\limits_{-\eta_{1}+\eta_{2}}^{\eta_{2}}}
%EndExpansion
6\,d\eta_{3}=6\eta_{1},\\%
%TCIMACRO{\dint \limits_{0}^{1-\eta_{1}}}%
%BeginExpansion
{\displaystyle\int\limits_{0}^{1-\eta_{1}}}
%EndExpansion
6\,d\eta_{3}=6\left(  1-\eta_{1}\right)  , &
%TCIMACRO{\dint \limits_{0}^{\eta_{2}}}%
%BeginExpansion
{\displaystyle\int\limits_{0}^{\eta_{2}}}
%EndExpansion
6\,d\eta_{3}=6\eta_{2}%
\end{array}
\]
and thus the joint density for $\eta_{1}$, $\eta_{2}$ is%
\[
\left\{
\begin{array}
[c]{lll}%
6\left(  1-\eta_{2}\right)   &  & \text{if }1-\eta_{2}<\eta_{1}<\eta_{2}\text{
and }1/2<\eta_{2}<1\text{,}\\
6\eta_{1} &  & \text{if }\eta_{1}<\eta_{2}<1-\eta_{1}\text{ and }0<\eta
_{1}<1/2\text{,}\\
6\left(  1-\eta_{1}\right)   &  & \text{if }1-\eta_{1}<\eta_{2}<\eta_{1}\text{
and }1/2<\eta_{1}<1\text{,}\\
6\eta_{2} &  & \text{if }\eta_{2}<\eta_{1}<1-\eta_{2}\text{ and }0<\eta
_{2}<1/2\text{.}%
\end{array}
\right.
\]
The sought-after joint density for $\alpha$, $\beta$ is therefore%
\[
\left\{
\begin{array}
[c]{lll}%
6\left(  \pi-\beta\right)  /\pi^{3} &  & \text{if }\pi-\beta<\alpha
<\beta\text{ and }\pi/2<\beta<\pi\text{,}\\
6\alpha/\pi^{3} &  & \text{if }\alpha<\beta<\pi-\alpha\text{ and }0<\alpha
<\pi/2\text{,}\\
6\left(  \pi-\alpha\right)  /\pi^{3} &  & \text{if }\pi-\alpha<\beta
<\alpha\text{ and }\pi/2<\alpha<\pi\text{,}\\
6\beta/\pi^{3} &  & \text{if }\beta<\alpha<\pi-\beta\text{ and }0<\beta<\pi/2
\end{array}
\right.
\]
and we call this the bivariate tent distribution (as opposed to pyramid
distribution, which already means something else \cite{Ke-quadrilat}). \ It is
clear that $\rho(\alpha,\beta)=0$ yet $\alpha$ and $\beta$ are dependent.%
%TCIMACRO{\FRAME{ftbpFU}{2.6401in}{2.66in}{0pt}{\Qcb{Density function for angle
%$\alpha$ in Section 2.}}{}{angles.eps}{\special{ language "Scientific Word";
%type "GRAPHIC";  maintain-aspect-ratio TRUE;  display "USEDEF";
%valid_file "F";  width 2.6401in;  height 2.66in;  depth 0pt;
%original-width 11.1041in;  original-height 11.188in;  cropleft "0";
%croptop "1";  cropright "1";  cropbottom "0";
%filename 'angles.eps';file-properties "XNPEU";}} }%
%BeginExpansion
\begin{figure}[ptb]%
\centering
\includegraphics[
height=2.66in,
width=2.6401in
]%
{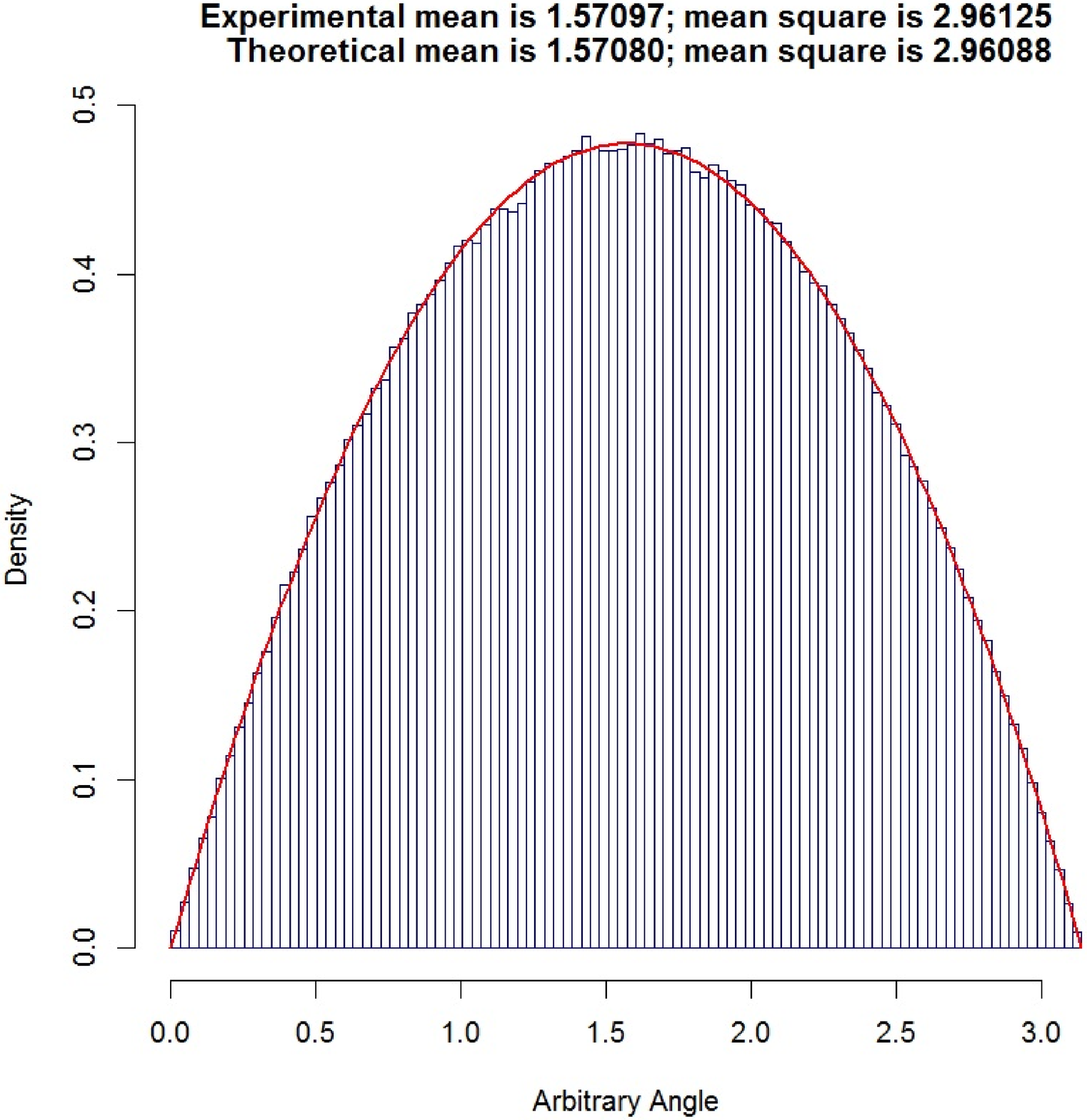}%
\caption{Density function for angle $\alpha$ in Section 2.}%
\end{figure}
%EndExpansion%
%TCIMACRO{\FRAME{ftbpFU}{4.9473in}{3.8904in}{0pt}{\Qcb{Density function for
%bivariate tent distribution on $[0,\pi]\times\lbrack0,\pi]$.}}{}%
%{tent.eps}{\special{ language "Scientific Word";  type "GRAPHIC";
%maintain-aspect-ratio TRUE;  display "USEDEF";  valid_file "F";
%width 4.9473in;  height 3.8904in;  depth 0pt;  original-width 13.255in;
%original-height 10.3985in;  cropleft "0";  croptop "1";  cropright "1";
%cropbottom "0";  filename 'tent.eps';file-properties "XNPEU";}} }%
%BeginExpansion
\begin{figure}[ptb]%
\centering
\includegraphics[
height=3.8904in,
width=4.9473in
]%
{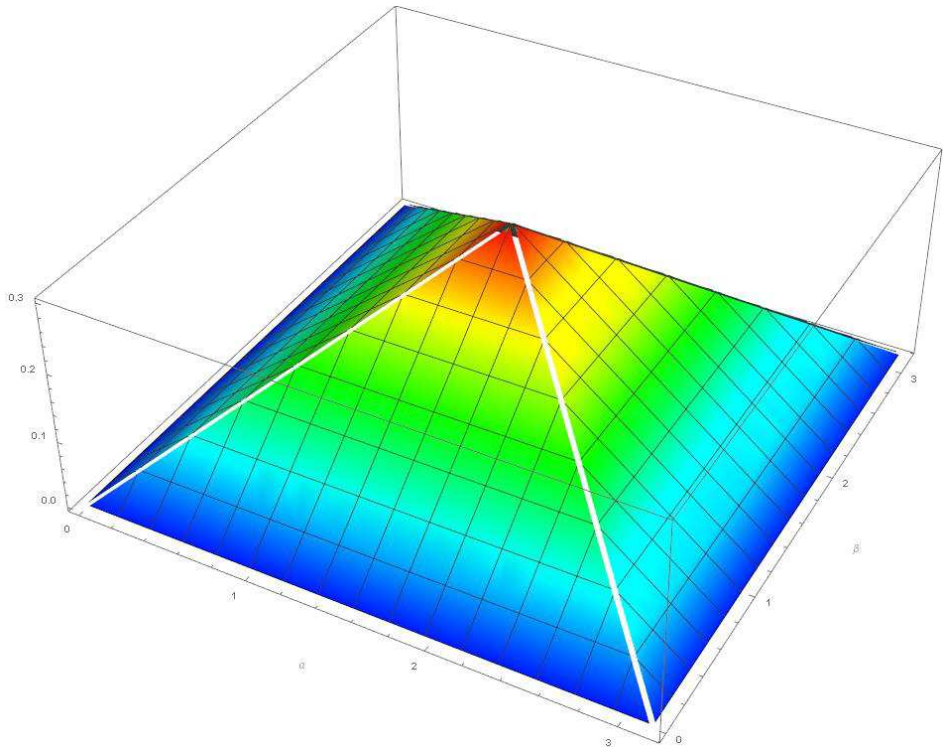}%
\caption{Density function for bivariate tent distribution on $[0,\pi
]\times\lbrack0,\pi]$.}%
\end{figure}
%EndExpansion

\section{Looking Back}

Given a uniform cyclic triangle, the joint density for two arbitrary angles
$\alpha$, $\beta$ is \cite{Ra-quadrilat, Mil-quadrilat, Mre-quadrilat}
\[
\left\{
\begin{array}
[c]{lll}%
2/\pi^{2} &  & \text{if }0<\alpha<\pi\text{, }0<\beta<\pi\text{ and }%
\alpha+\beta<\pi,\\
0 &  & \text{otherwise}%
\end{array}
\right.
\]
and trivially $\rho(\alpha,\beta)=-1/2$. \ Let $\Delta$ denote the isosceles
triangular support of this distribution. \ Let $a$ denote the side opposite
$\alpha$ and $b$ denote the side opposite $\beta$. $\ $From%
\[
\left(
\begin{array}
[c]{c}%
\alpha\\
\beta
\end{array}
\right)  \mapsto\left(
\begin{array}
[c]{c}%
2\sin(\alpha)\\
2\sin(\beta)
\end{array}
\right)  =\left(
\begin{array}
[c]{c}%
a\\
b
\end{array}
\right)  ,
\]
we have Jacobian determinant $4\cos(\alpha)\cos(\beta)$ and preimages%
\[%
\begin{array}
[c]{ccc}%
\left(
\begin{array}
[c]{c}%
\arcsin\left(  \tfrac{a}{2}\right) \\
\arcsin\left(  \tfrac{b}{2}\right)
\end{array}
\right)  , &  & \left(
\begin{array}
[c]{c}%
\pi-\arcsin\left(  \tfrac{a}{2}\right) \\
\arcsin\left(  \tfrac{b}{2}\right)
\end{array}
\right)
\end{array}
\]
if $b<a$ and%
\[%
\begin{array}
[c]{ccc}%
\left(
\begin{array}
[c]{c}%
\arcsin\left(  \tfrac{a}{2}\right) \\
\arcsin\left(  \tfrac{b}{2}\right)
\end{array}
\right)  , &  & \left(
\begin{array}
[c]{c}%
\arcsin\left(  \tfrac{a}{2}\right) \\
\pi-\arcsin\left(  \tfrac{b}{2}\right)
\end{array}
\right)
\end{array}
\]
if $a<b$. \ Reason: if $b<a$, then $\arcsin(b/2)<\arcsin(a/2)$ and hence both
preimages fall in $\Delta$ because $\left[  \pi-\arcsin(a/2)\right]
+\arcsin(b/2)<\pi$. \ No other preimages exist when $b<a$ because
$\arcsin(a/2)+\left[  \pi-\arcsin(b/2)\right]  >\pi$ and $\left[  \pi
-\arcsin(a/2)\right]  +\left[  \pi-\arcsin(b/2)\right]  >\pi$. \ Likewise for
$a<b$. \ 

The joint density for $a$ and $b$ is thus
\[
\left\{
\begin{array}
[c]{lll}%
\dfrac{4}{\pi^{2}}\dfrac{1}{\sqrt{4-a^{2}}}\dfrac{1}{\sqrt{4-b^{2}}} &  &
\text{if }0<a<2\text{ and }0<b<2\text{,}\\
0 &  & \text{otherwise}%
\end{array}
\right.
\]
which implies that sides $a$, $b$ are independent even though they are related
so easily (via the sine function) to the dependent angles $\alpha$, $\beta$.
\ As far as is known, this observation is new. \ We mention that the remaining
side $c$ satisfies%
\[
c=\left\{
\begin{array}
[c]{lll}%
\frac{1}{2}\left(  a\sqrt{4-b^{2}}+b\sqrt{4-a^{2}}\right)  &  & \text{with
probability }1/2,\\
\frac{1}{2}\left\vert a\sqrt{4-b^{2}}-b\sqrt{4-a^{2}}\right\vert  &  &
\text{with probability }1/2
\end{array}
\right.
\]
for completeness' sake.

\section{Looking Forward}

The polygonal angles $\alpha$, $\beta$, $\gamma$, $\delta$ associated with a
uniform cyclic $5$-gon can be studied via the Dirichlet$(1,1,1,1;1)$
distribution on a $4$-dimensional simplex \cite{Ra-quadrilat, PC-quadrilat,
NT-quadrilat}:%
\[
\left\{
\begin{array}
[c]{lll}%
24 &  & \text{if }0<\xi_{1}<1,\text{ }0<\xi_{2}<1,\text{ }0<\xi_{3}<1,\text{
}0<\xi_{4}<1\text{ and }\xi_{1}+\xi_{2}+\xi_{3}+\xi_{4}<1,\\
0 &  & \text{otherwise}%
\end{array}
\right.
\]
and calculation of the joint density for $\eta_{1}=\xi_{1}+\xi_{2}+\xi_{3}$,
$\eta_{2}=\xi_{2}+\xi_{3}+\xi_{4}$, $\eta_{3}=1-\xi_{1}-\xi_{2} $, $\eta
_{4}=1-\xi_{2}-\xi_{3}$. \ Omitting elaborate details, we obtain the density
to be $24$ when%
\[
\max\{1-\eta_{1},1-\eta_{2}\}<\eta_{4}<\min\{2-\eta_{1}-\eta_{2},2-\eta
_{1}-\eta_{3}\}\text{ and }1<\eta_{1}+\eta_{3}<2
\]
and $0$ otherwise. \ It follows that%
\[%
\begin{array}
[c]{ccccccc}%
\rho(\alpha,\beta)=1/6, &  & \rho(\alpha,\gamma)=-2/3, &  & \rho(\alpha
,\delta)=-2/3, &  & \rho(\alpha,\varphi)=1/6
\end{array}
\]
where $\varphi=3\pi-\alpha-\beta-\gamma-\delta$. \ In particular, adjacent
angles are positively correlated and non-adjacent angles are negatively correlated.

For a uniform cyclic $6$-gon, we conjecture that%
\[%
\begin{array}
[c]{lllll}%
\rho(\alpha,\beta)=1/4, &  & \rho(\alpha,\gamma)=-1/2, &  & \rho(\alpha
,\delta)=-1/2,\\
\rho(\alpha,\varphi)=-1/2, &  & \rho(\alpha,\psi)=1/4 &  &
\end{array}
\]
where $\varphi=2\pi-\alpha-\gamma$ and $\psi=2\pi-\beta-\delta$. \ Again,
adjacent angles are positively correlated and non-adjacent angles are
negatively correlated. \ The fact that $\delta$ is opposite $\alpha$ seems not
to affect its correlation with $\alpha$, relative to either $\gamma$ or
$\varphi$.

\section{Area}

Given a uniform cyclic triangle, moments of area $2\sin(\alpha)\sin(\beta
)\sin(\alpha+\beta)$ are computed by use of the joint angle density:
\[
\frac{2}{\pi^{2}}%
%TCIMACRO{\dint \limits_{0}^{\pi}}%
%BeginExpansion
{\displaystyle\int\limits_{0}^{\pi}}
%EndExpansion%
%TCIMACRO{\dint \limits_{0}^{\pi-\beta}}%
%BeginExpansion
{\displaystyle\int\limits_{0}^{\pi-\beta}}
%EndExpansion
2\sin(\alpha)\sin(\beta)\sin(\alpha+\beta)\,d\alpha\,d\beta=\frac{3}{2\pi},
\]%
\[
\frac{2}{\pi^{2}}%
%TCIMACRO{\dint \limits_{0}^{\pi}}%
%BeginExpansion
{\displaystyle\int\limits_{0}^{\pi}}
%EndExpansion%
%TCIMACRO{\dint \limits_{0}^{\pi-\beta}}%
%BeginExpansion
{\displaystyle\int\limits_{0}^{\pi-\beta}}
%EndExpansion
4\sin^{2}(\alpha)\sin^{2}(\beta)\sin^{2}(\alpha+\beta)\,d\alpha\,d\beta
=\frac{3}{8}.
\]
The density for area itself is $8xK\left(  4x^{2}\right)  $, where%
\begin{align*}
K(y)  &  =\frac{1}{4\pi^{3}}\frac{1}{\sqrt{y}}\left\{  \Gamma\left(  \frac
{1}{3}\right)  ^{3}\left(  \frac{4y}{27}\right)  ^{-1/6}\,_{2}F_{1}\left(
\dfrac{1}{3},\dfrac{1}{3},\dfrac{2}{3},\frac{4y}{27}\right)  -\right. \\
&  \ \ \ \left.  3\Gamma\left(  \frac{2}{3}\right)  ^{3}\left(  \frac{4y}%
{27}\right)  ^{1/6}\,_{2}F_{1}\left(  \dfrac{2}{3},\dfrac{2}{3},\dfrac{4}%
{3},\frac{4y}{27}\right)  \right\}  ,
\end{align*}
$_{2}F_{1}$ is the Gauss hypergeometric function and $0<y<27/4$. \ This
formula corrects that which appears in Case III of \cite{MT-quadrilat}.

Given a uniform cyclic quadrilateral, we conjecture that the joint density for
angles $\alpha$, $\beta$, $\omega$ is%
\[
f(\alpha,\beta,\omega)=\left\{
\begin{array}
[c]{lll}%
3/\pi^{3} &  & \text{if }\alpha+\beta>\omega\text{, }\alpha+\omega
>\beta\text{, }\beta+\omega>\alpha\text{ and }\alpha+\beta+\omega<2\pi,\\
0 &  & \text{otherwise.}%
\end{array}
\right.
\]
It can be shown that, assuming the formula for $f$ is valid, any two angles
from the list $\alpha$, $\beta$, $\omega$ are distributed according to the
bivariate tent density. \ Our conjecture is consistent with computer
simulation, but a rigorous proof is open. \ From this, we obtain area moments%
\[%
%TCIMACRO{\dint \limits_{0}^{\pi}}%
%BeginExpansion
{\displaystyle\int\limits_{0}^{\pi}}
%EndExpansion%
%TCIMACRO{\dint \limits_{0}^{\pi}}%
%BeginExpansion
{\displaystyle\int\limits_{0}^{\pi}}
%EndExpansion%
%TCIMACRO{\dint \limits_{0}^{\pi}}%
%BeginExpansion
{\displaystyle\int\limits_{0}^{\pi}}
%EndExpansion
2\sin(\alpha)\sin(\beta)\sin(\omega)f(\alpha,\beta,\omega)\,d\alpha
\,d\beta\,d\omega=\frac{3}{\pi},
\]%
\[%
%TCIMACRO{\dint \limits_{0}^{\pi}}%
%BeginExpansion
{\displaystyle\int\limits_{0}^{\pi}}
%EndExpansion%
%TCIMACRO{\dint \limits_{0}^{\pi}}%
%BeginExpansion
{\displaystyle\int\limits_{0}^{\pi}}
%EndExpansion%
%TCIMACRO{\dint \limits_{0}^{\pi}}%
%BeginExpansion
{\displaystyle\int\limits_{0}^{\pi}}
%EndExpansion
4\sin^{2}(\alpha)\sin^{2}(\beta)\sin^{2}(\omega)f(\alpha,\beta,\omega
)\,d\alpha\,d\beta\,d\omega=\frac{1}{2}+\frac{105}{16\pi^{2}}
\]
which again is consistent with experiment. \ The mean area for quadrilaterals
is twice that for triangles. No formula for the density of area itself is
known.%
%TCIMACRO{\FRAME{ftbpFU}{4.3742in}{4.7565in}{0pt}{\Qcb{Tetrahedral support for
%$f$, with vertices $(0,0,0)$, $(0,\pi,\pi)$, $(\pi,0,\pi)$, $(\pi,\pi,0)$.}}%
%{}{tetra0.eps}{\special{ language "Scientific Word";  type "GRAPHIC";
%maintain-aspect-ratio TRUE;  display "USEDEF";  valid_file "F";
%width 4.3742in;  height 4.7565in;  depth 0pt;  original-width 9.692in;
%original-height 10.5447in;  cropleft "0";  croptop "1";  cropright "1";
%cropbottom "0";  filename 'tetra0.eps';file-properties "XNPEU";}} }%
%BeginExpansion
\begin{figure}[ptb]%
\centering
\includegraphics[
height=4.7565in,
width=4.3742in
]%
{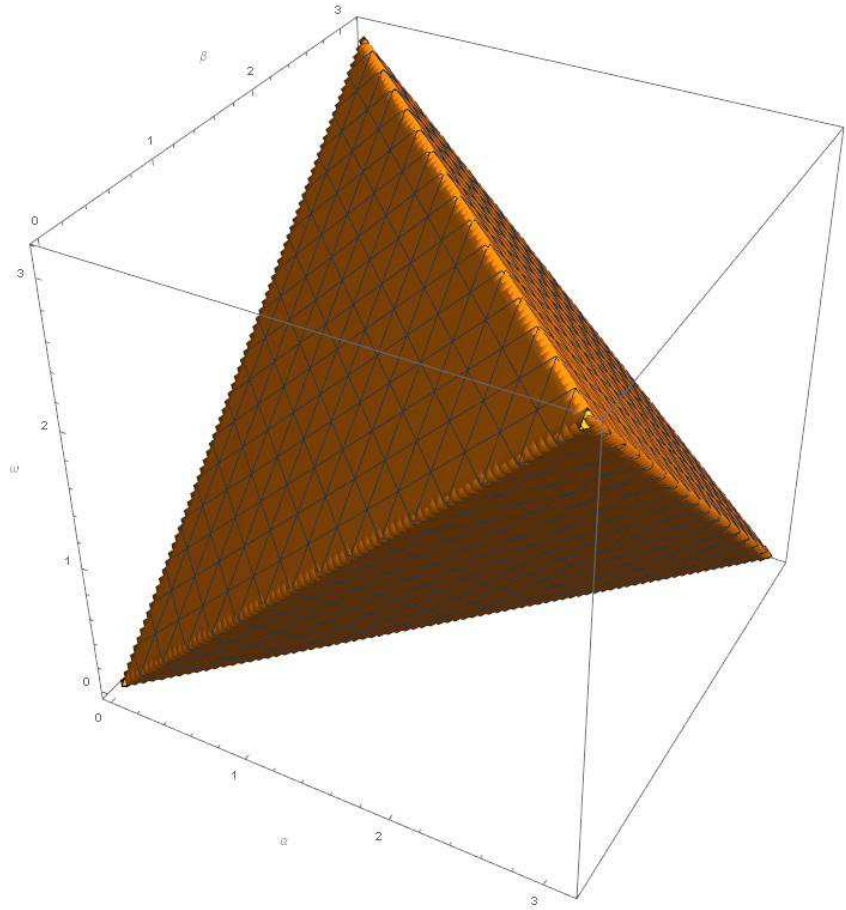}%
\caption{Tetrahedral support for $f$, with vertices $(0,0,0)$, $(0,\pi,\pi)$,
$(\pi,0,\pi)$, $(\pi,\pi,0)$.}%
\end{figure}
%EndExpansion

The problem with angles is that we do not know a suitable way of relating
$\omega$ with parameters $\theta_{1}$, $\theta_{2}$, $\theta_{3}$, $\theta
_{4}$. \ For a cyclic quadrilateral with successive sides $a$, $b$, $c $, $d$,
formulas like \cite{DR-quadrilat}
\[%
\begin{array}
[c]{ccc}%
\tan\left(  \dfrac{\alpha}{2}\right)  =\sqrt{\dfrac{(-a+b+c+d)(a-b+c+d)}%
{(a+b-c+d)(a+b+c-d)}} &  & \text{where }\alpha\text{ is angle between }a\text{
and }b\text{,}%
\end{array}
\]%
\[%
\begin{array}
[c]{ccc}%
\tan\left(  \dfrac{\beta}{2}\right)  =\sqrt{\dfrac{(a-b+c+d)(a+b-c+d)}%
{(-a+b+c+d)(a+b+c-d)}} &  & \text{where }\beta\text{ is angle between }b\text{
and }c\text{,}%
\end{array}
\]%
\[%
\begin{array}
[c]{ccc}%
\tan\left(  \dfrac{\omega}{2}\right)  =\sqrt{\dfrac{(a-b+c+d)(a+b+c-d)}%
{(-a+b+c+d)(a+b-c+d)}} &  & \text{where }\omega\text{ is angle between
diagonals}%
\end{array}
\]
suggest an alternative approach to solution, but the path seems very complicated.

\section{Acknowledgements}

I am indebted to Chi Zhang for her hand calculations in Sections 2 and 4
(specifically, those involving $\xi$s and $\eta$s). I am also grateful to
Guo-Liang Tian, Serge Provost and Paul Kettler for helpful discussions.

\end{document}